\documentclass[a4paper,11pt,dvips]{article}

\usepackage{amsfonts}
\usepackage{amsmath}
\usepackage{amsthm}
\usepackage[auth-sc,affil-sl]{authblk}
\usepackage{bm}
\usepackage{multirow}
\usepackage{verbatim}
\usepackage{enumitem}

\renewcommand
\baselinestretch{1.0} 
\usepackage{url}
\usepackage{natbib}
\bibpunct{(}{)}{,}{a}{}{,}
\usepackage{graphicx}
\usepackage{lscape}
\usepackage[usenames]{color}
\usepackage{sansmath}
\usepackage{subfigure}
\usepackage{rotating}

\DeclareMathAlphabet{\mathpzc}{OT1}{pzc}{m}{it}

\def\V{\mathrm{Var}}

\newcommand\independent{\protect\mathpalette{\protect\independenT}{\perp}}
\def\independenT#1#2{\mathrel{\rlap{$#1#2$}\mkern2mu{#1#2}}}

\def\Expect{\mathbb{{E}}}
\def\P{\text{Pr}}
\newtheorem{Thm}{Theorem}

\newcommand \address[1]{\gdef \@address{#1}}
\makeatletter
\long\def\@footnotetext#1{\insert\footins{\def\baselinestretch{1.2}\footnotesize
\interlinepenalty\interfootnotelinepenalty
\splittopskip\footnotesep \splitmaxdepth \dp\strutbox
\floatingpenalty \@MM \hsize\columnwidth \@parboxrestore
\edef\@currentlabel{\csname
p@footnote\endcsname\@thefnmark}\@makefntext
{\rule{\z@}{\footnotesep}\ignorespaces #1\strut}}}

\long\def\symbolfootnote[#1]#2{\begingroup%
\def\thefootnote{\fnsymbol{footnote}}\footnote[#1]{#2}\endgroup}

\providecommand{\keywords}[1]{\textit{Keywords:} #1}

\def\maketitle{%
  \null
  \thispagestyle{empty}%
  \begin{center}\leavevmode
    \normalfont
    {\LARGE \bf \@title\par}%
    {\normalsize \@author\par}%
    \vskip 0.05 cm
    \vskip 0.05cm
    {\normalsize \@date\par}%
  \end{center}%
}

\makeatletter
\newcommand{\institute}[1]{\newcommand{\@institute}{#1}}

\renewcommand\textsl{\textcolor{blue}}

\begin{document}
\title{Doubly robust dose-response estimation for continuous treatments via generalized propensity score augmented outcome regression}
\author[1]{Daniel J. Graham}
\affil[1]{Corresponding author: Department of Civil Engineering, Imperial College London, London, SW7 2AZ, UK. Email: \texttt{d.j.graham@imperial.ac.uk}} 
%\affil[ ]{Email: \texttt{d.j.graham@imperial.ac.uk}}
%
\author[2]{Emma J. McCoy}
\affil[2]{Department of Mathematics, Imperial College London, London, UK}

\author[3]{David A. Stephens}
\affil[3]{Department of Mathematics and Statistics, McGill University, Montreal, Canada}

\date{}

\maketitle
\begin{abstract}
This paper constructs a doubly robust estimator for continuous dose-response estimation. An outcome regression model is augmented with a set of inverse generalized propensity score covariates to correct for potential misspecification bias. From the augmented model we can obtain consistent estimates of mean average potential outcomes for distinct strata of the treatment. A polynomial regression is then fitted to these point estimates to derive a Taylor approximation to the continuous dose-response function. The bootstrap is used for variance estimation. Analytical results and simulations show that our approach can provide a good approximation to linear or nonlinear dose-response functions under various sources of misspecification of the outcome regression or propensity score models. Efficiency in finite samples is good relative to minimum variance consistent estimators.
\end{abstract}
.

\section{Introduction}

The typical set up in causal inference problems is one in which the data available for estimation are realisations of a random vector, $Z_i = (Y_i,D_i,X_i)$, $i=1,...,n$, where for the $i$-th unit of observation $Y_i$ denotes a response, $D_i$ the treatment received, and $X_i$ a vector of pre-treatment covariates. The objective is to estimate an average treatment effect (ATE), or in other words, the difference in response that would occur under different treatment assignments averaged over the population. However, since the treatment is usually not assigned randomly, simple comparisons of mean responses across different treatment groups will not in general reveal a `causal' effect due to confounding.

Consistent causal estimates of ATEs can, however, be obtained if the vector of covariates $X_i$ is sufficient to ensure conditional independence of potential outcomes and treatment assignment. Using the notation of \citet{Tsiatis/Davidian:2007}, we write joint densities of the observed data in the form
\[f_Z(z)=f_{Y|D,X}(y_i|d_i,x_i)f_{D|X}(d_i|x_i)f_X(x_i).\]
Well known causal ATE estimators can be derived via an outcome regression (OR) model for $\Expect(Y_i|D_i,X_i)$, the mean of the conditional density of the response given the covariates and treatment status; or via a propensity score (PS) model, $\pi(D_i|X_i)$, for the treatment assignment (or exposure) mechanism.

\textit{Doubly robust} (DR) estimators combine the OR and PS models, usually by weighting or augmenting the OR model with covariates derived from the inverse propensity score, to obtain estimates of ATEs that are consistent and asymptotically normal when either the OR or the PS model are correctly specified. DR estimators been studied extensively in the context of binary treatments \citep[e.g.][]{Robins:2000,Robins/et/al:2000c,Robins/Rotnitzky:2001,VanDerLaan/Robins:2003,Lunceford/Davidian:2004,Bang/Robins:2005,Kang/Schafer:2007}.

In this paper we extend the binary DR augmented regression approach of \citet{Scharfstein/et/al:1999} to derive a Taylor approximation for continuous dose-response functions. We present theoretical results which explain the DR properties of our approach and simulation results which show that the estimator can be used to recover a good approximation to linear or nonlinear dose-response functions.

The paper is structured as follows. Section two explains some general principles of PS based continuous dose-response estimation. Section three presents our DR model for continuous treatments and demonstrates the consistency properties of this estimator. Simulation results are presented in section four. Conclusions are drawn in the final section.

\section{Potential outcomes and dose-response estimation for continuous treatments}

Throughout the paper we use the convention of upper-case to denote a random variable and lower case for its observed value. In the single time-point setting, let $d \in \mathcal{D}  \subseteq \mathbb{R}$, denote a given value (dose) of treatment and let $d_i$ denote the value of the treatment actually assigned to unit $i$.  For each unit $i$ let there be a potential outcome $Y_i(d)$ associated with each dose, such that $\mathcal{Y}_i=\{Y_i(d): d \in \mathcal{D}\}$ denotes the full set of potential outcomes. In the data we observe random variables describing the actual dose received $D_i$, the outcome associated with that dose $Y_i(D_i)$, and a set of pre-treatment covariates $X_i$.  We define the indicator for receipt of dose $d$ as
\begin{equation*}
I_{d}(D_i)\stackrel{}{=} \left\{
\begin{array}{ll} 1 & \text{if  } D_i = d \\
0 & \text{if  } D_i \neq d
\end{array} \right.,
\end{equation*}
and similarly the indicator for receipt of the dose actually received is $I_{d_i}(D_i)$. We assume that the relationship between actual and potential outcomes obeys the stable unit treatment value assumption (SUTVA) such that $Y_i \equiv I_{d}(D_i)Y_i(d)$ for all $d \in \mathcal{D}$, for all $Y_i (d) \in \mathcal{Y}_i$, and for $i=1,...,n$.

Our starting point for the analysis of continuous treatments within a potential outcomes framework is the assumption of weak conditional independence introduced by \citet{Imbens:2000}, which requires that $Y_i(d) \independent I_{d}(D_i)|X_i$ for all $d \in \mathcal{D}$. If weak conditional independence holds we can obtain causal estimates from the observed data without knowledge of the full range of potential outcomes because
\begin{equation}\label{APO}
\Expect \left[Y_i(d)|X_i \right]=\Expect \left[Y_i(d)|I_{d}(D_i),X_i \right]=\Expect \left[Y_i|I_{d}(D_i),X_i \right],
\end{equation}
and taking the expectation over covariates $X_i$ gives the Average Potential Outcome (APO) at dose $d$, $\mu(d)=\Expect \left[Y_i(d) \right]$, commonly referred to as the dose-response function.

An OR based approach could be used to estimate (\ref{APO}) by specifying a mean response model $\Expect(Y_i|D_i,X_i)=\Psi^{-1} \{m(D_i, X_i; \beta) \}$, for some regression function $m()$ with known link function $\Psi$ and unknown parameter vector $\beta$. If the assumed model provides a correct specification of the true mean response then in principle a consistent estimate of the dose-response at dose $d$ can then be obtained using
\begin{equation} \label{DROR}
\widehat{\mu}_{OR}(d)=n^{-1} \sum_{i=1}^{n} \Psi^{-1} \left\{m(d, X_i; \widehat{\beta})\right\},
\end{equation}
that is, the average across all data of the predicted values evaluated at $D_i=d$

As in the case of binary treatments, however, estimation of continuous dose-response functions is often made more tractable by working with a scalar PS rather than the potentially high dimensional covariate vector. For continuous treatments, \citet{Imbens:2000} and \citet{Hirano/Imbens:2004} introduce the generalized propensity score (GPS). Let $r(d,x_i)=f_{D|X}(d|x_i)$ denote the conditional density function for receiving a particular dose of the treatment given pre-treatment variables $X_i=x_i$. The observed GPS, which we denote $R_i=r(D_i,X_i)$, is a random variable comprising values from this conditional density evaluated at the level of treatment actually received (i.e. $D_i$) given $X_i=x$. In addition, we also define the family of random variables indexed by $d$, $R_{d,i}=r(d,X_i)$, as values from the conditional density of receiving a particular dose $d$ given $X_i$. Clearly when $D_i=d$, $R_i =R_{d,i}$.

\citet{Hirano/Imbens:2004} show that the GPS provides a bias removal strategy in the context of continuous treatments such that
\begin{equation} \label{DRGPS}
\mu(d)\equiv \Expect[Y_i(d)]\equiv \Expect_{X}[\Expect(Y_i(d)|X)]=\Expect_{X}[\Expect[Y_i(d)|r(d,X_i)]]=\Expect_{X}[\Expect[Y_i|I_{d}(D_i),r(d,X_i)]].
\end{equation}

Dose response estimation can therefore proceed via OR or GPS based models. The consistency of the estimator (\ref{DROR}) relies on correct specification of an OR model while that of (\ref{DRGPS}) relies on correct specification of the GPS. In the next section, we develop a DR estimator for continuous and potentially nonlinear dose-response functions by extending the DR approach for binary treatments of \citet{Scharfstein/et/al:1999}. 

\section{Doubly robust estimation for continuous treatments}

Dose-response estimation is challenging because the relationship may be nonlinear and the bias induced by sources of misspecification non-constant across doses. In constructing an appropriate dose-response estimator, flexibility is required to address the inherent trade off between fidelity to these nonlinearities and consistency of the estimated causal quantities.

Our proposed approach involves two main steps. First we specify an augmented OR model to obtain consistent DR estimates of mean APOs for strata of the treatment. We then use polynomial regression to recover point estimates of the dose-response curve with standard errors calculated via the bootstrap.

\subsection{Estimating equations for DR mean APO estimation}

In the first step, our estimand of interest is the mean APO for strata of the treatment.  We proceed as follows:

\begin{enumerate}[label=({\roman*}),topsep=-2.4em,itemsep=-0.5em,parsep=1em,partopsep=2em]
\item We define $Q$, $q=(1,...,Q)$, strata over the range of $d$ and use $\mathcal{D}_q$ to denote treatment stratum $q$ $(\mathcal{D}_q \subset \mathcal{D} \subseteq \mathbb{R}$). We index distinct treatment levels within each stratum by $j=(1,...,J)$, such that individual elements of $\mathcal{D}_q$ are denoted $d_{qj}$. Our objective is to estimate
\begin{equation}
\mu\left(\mathcal{D}_q\right)=\Expect_j\left[\Expect\left\{Y_i(d_{qj})\right\}\right].
\end{equation}

\item Given some assumed density, $f_{D|X}(d|x_i, \alpha)$, with parameters vector $\alpha$, and for some small $\delta$ defined around the dose of interest, we define a probability representation of the GPS as
\begin{align}\label{PGPS}
\pi(d|X_i;\alpha)&=\P(d-\delta \leq D_i \leq d+\delta|X_i,\alpha)\\ \nonumber
&=\int_{d-\delta}^{d+\delta}f_{D|X}(t|x_i,\alpha)dt \\ \nonumber
&= \Expect\left[I_{d}(D_i)|X_i=x\right],
\end{align}
which we refer to as the \textit{probabilistic generalized propensity score} (PGPS). 

\item Specifying an appropriate regression model for $f_{D|X}(d|x_i, \alpha)$, we use observed doses $d_i$ and observed covariates $x_i$ to obtain consistent estimates of $\widehat{\alpha}$ and calculate observed PGPSs: $\widehat{\pi}(D_i|X_i;\widehat{\alpha})$. 

\item Denoting membership of treatment stratum $q$ using the indicator $I_q(D_i)$, we specify an augmented OR (AOR) model as
\begin{equation} \label{DRM}
e\left\{D_i,X_i,\widehat{\kappa}_i(D_i,X_i);\xi\right\}=\Psi^{-1} \left\{m\left(D_i,X_i; \beta \right)+\sum_{q=1}^Q \frac{\varphi_{q} I_q(D_i)}{\widehat{\pi}(D_i|X_i;\widehat{\alpha})}\right\}=\Psi^{-1} \left\{m\left(D_i,X_i,\widehat{\kappa}_i(D_i,X_i);\xi \right) \right\}
\end{equation}
say, where $\varphi=(\varphi_1,...,\varphi_Q)$ is a $Q$ dimensional parameter vector for the inverse PGPS covariates, $\xi=(\beta,\varphi)$, and $\widehat{\kappa}_i(D_i,X_i)$ is an ($n \times Q$) matrix each column of which contains a covariate for treatment stratum $q$
\[\frac{I_q(D_i)}{\widehat{\pi}(D_i|X_i;\widehat{\alpha})}.\] 
We obtain estimates of $\xi=(\beta,\varphi)$ as solutions to estimating equations of the form
\begin{equation}
\sum_{i=1}^{n}\frac{1}{\phi}\frac{\partial e\left\{d_i,x_i,\widehat{\kappa}_i(d_i,x_i); \xi\right\}}{\partial \xi^{\sf{T}}}\left[y_i-e\left\{d_i,x_i,\widehat{\kappa}_i(d_i,x_i; \xi)\right\} \right]=0
\end{equation}
where $\phi$ is a working variance matrix for $\V[Y_i|D_i,X_i,\kappa_i(D_i,X_i)]$. 

\item For each stratum, we then calculate the means of the predicted values of the AOR model evaluated at each level of $\mathcal{D}_q$. Our estimator is the average of these means
\begin{equation} \label{MUDR1}
\widehat{\mu}_{DR}\left(\mathcal{D}_q\right)=\frac{1}{J}\sum_{j=1}^J \left[\frac{1}{n}\sum_{i=1}^{n} \Psi^{-1} \left\{m\left(d_{qj},X_i; \widehat{\beta} \right) + \frac{\varphi_{q}} {\widehat{\pi}(d_{qj}|X_i; \widehat{\alpha})}\right\} \right].
\end{equation}
\item Finally, we slide the partition to define new strata and repeat steps (i) to (v) thus increasing coverage of mean APOs on the dose-response function.
\end{enumerate}

For valid inference, the calculations must be based on a defined region of dose, $\mathcal{C}\subseteq \mathcal{D}$, within which there is common support by treatment status in the covariate distributions. The common support requirement is met when for any subset of $\mathcal{C}$, say $\mathcal{A}\subseteq \mathcal{C}$, $\P(d \in \mathcal{A}|X_i=x)>0$ for all $d$ and $x$ and $\mathcal{A}\subseteq \mathcal{C}$. 

The estimator is doubly robust in the sense that it yields a consistent estimate of the mean APO in each stratum if one of two parametric restrictions hold: 1) the OR model $\Psi^{-1}\{m(X_i, D_i; \widehat{\xi})\}$ provides an consistent estimate of the true conditional expectation $\Expect(Y_i|D_i,X_i)$, or; 2) $\widehat{\pi}(d|X_i;\widehat{\alpha})$ provides a consistent estimator of the true PGPS $\Expect\left[I_{d}(D_i)|X_i=x\right]$. These restrictions, and their implications for mean APO estimation, are described in the following two theorems on redundant conditioning and bias correction.

\begin{Thm}(Effects of redundant conditioning on APO estimation). If the OR model $\Psi^{-1}\{m(X_i, D_i; \widehat{\beta})\}$ provides a consistent estimate of $\Expect[Y_i|D_i,X_i]$, then the mean of the predicted values from this model evaluated at $D_i=d$ will provide a consistent estimate of $\Expect[Y_i(d)]$ provided conditional independence between treatment status and response holds. The mean predicted values from the AOR model, evaluated at $D_i=d$, will also produce a consistent estimate of $\Expect[Y_i(d)]$ under conditional independence, but the AOR model will be less efficient.    \\
\\
{\bf Proof.}
The true conditional expectation to be estimated is $\Expect[Y_i|D_i,X_i]$, but instead we estimate a model for $\Expect\left[Y_i|D_i,X_i,\widehat{\kappa}_i(D_i,X_i)\right]$. If we have conditional independence given covariates $X_i$, and if the SUTVA holds, then
\[\Expect[Y_i(d)|X_i,\widehat{\kappa}_i(D_i,X_i)]=\Expect[Y_i|d,X_i,\widehat{\kappa}_i(D_i,X_i)],\]
regardless of redundant conditioning on $\widehat{\kappa}_i(D_i,X_i)$. The dose-response can then be obtained by taking expectations
\[\Expect[Y_i(d)]=\Expect_{X_i, \widehat{\kappa}_i(D_i,X_i)}[\Expect[Y_i|d,X_i,\widehat{\kappa}_i(D_i,X_i)]]=\Expect_{X_i}[\Expect[Y_i|d,X_i]].\]
However, the AOR model will have higher variance than the OR model because
\begin{align*}
\V\left(\Expect\left[Y_i(d)|X_i,\widehat{\kappa}_i(D_i,X_i)\right]\right)=&\V\left(\Expect\left[\Expect\left[Y_i(d)|X_i,\widehat{\kappa}_i(D_i,X_i)\right]|\widehat{\kappa}_i(D_i,X_i)\right] \right)\\
&+\Expect\left[\V \left(\Expect\left[Y_i(d)|X_i,\widehat{\kappa}_i(D_i,X_i)\right]|\widehat{\kappa}_i(D_i,X_i) \right) \right]\\
=&\V\left(\Expect\left[Y_i(d)|X_i\right] \right)+\Expect\left[\V \left(\Expect\left[Y_i(d)|X_i,\widehat{\kappa}_i(D_i,X_i)\right]|\widehat{\kappa}_i(D_i,X_i) \right) \right]
\end{align*}
and since $\Expect\left[\V \left(\Expect\left[Y_i(d)|X_i,\widehat{\kappa}_i(D_i,X_i)\right]|\widehat{\kappa}_i(D_i,X_i) \right) \right] \geq  0$, then $\V\left(\Expect[Y_i(d)|X_i,\widehat{\kappa}_i(D_i,X_i))\right] \geq  \V\left(\Expect[Y_i(d)|X_i]\right)$. 
\end{Thm}

\begin{Thm}(The bias correction property). If the estimated propensity score $\widehat{\pi}(d|X_i;\widehat{\alpha})$ provides a consistent estimate of the true propensity score $\pi(d|X_i;\widehat{\alpha})=\Expect\left[I_{d}(D_i)|X_i=x\right]$, then the mean predicted values from the AOR model, averaged within some stratum of the treatment, will provide consistent estimates of the mean APO for that stratum, $\mu\left(\mathcal{D}_q\right)$, under misspecification of the OR model.
\\
\\
{\bf Proof.}
Parameter estimates $\widehat{\varphi}$ are obtained from the AOR model as solution to estimating equations of the form
\[n^{-1}\sum_{i=1}^{n} \frac{I_q(D_i)}{\widehat{\pi}_i(D_i|X_i;\widehat{\alpha})} \cdot \left[y_i-\Psi^{-1} \left\{m\left(d_i,x_i; \beta \right) +\sum_{q=1}^Q \varphi_{q} \frac{I_q(d_i)} {\widehat{\pi}_i(d_i|x_i;\widehat{\alpha})}\right\}\right]=0,\]
which converge in probability to
\[\Expect\left(\frac{1}{\widehat{\pi}_i(D_i|X_i;\widehat{\alpha})}\cdot\left[Y_i-\Psi^{-1} \left\{m\left(D_i,X_i; \beta\right)+\sum_{q=1}^Q \varphi_{q} \frac{I_q(D_i)} {\widehat{\pi}_i(D_i|X_i;\widehat{\alpha})}\right\}\right]\biggr|I_q(D_i)=1\right).\]
Writing the AOR model as $\Psi^{-1} \left\{\cdot \right\}$ and taking expectations over $X_i$ and $D_i$ gives
\[ \Expect_{X_i,D_i}\left[\Expect\left(I_{d_i}(D_i)\cdot \frac{1}{\widehat{\pi}_i(D_i|X_i;\widehat{\alpha})} \cdot\left[Y_i-\Psi^{-1} \left\{\cdot \right\}\right]\biggr| I_q(D_i)=1,X_i,D_i\right)\right].\]
From weak conditional independence and the SUTVA
\[\Expect_{X_i,D_i}\left[\Expect\left[I_{d_i}(D_i)|X_i\right]\cdot I_q(D_i)\cdot \Expect\left(\frac{1}{\widehat{\pi}_i(D_i|X_i;\widehat{\alpha})} \left[Y_i(D_i)-\Psi^{-1} \left\{\cdot\right\}\right]\biggr|I_q(D_i)=1,X_i\right)\right],\]
and since $\widehat{\pi}_i$ is a function of $X_i$ we use the pull through property to give
\[\Expect_{X_i,D_i}\left[\Expect\left[I_{d_i}(D_i)|X_i\right] \cdot \frac{1}{\widehat{\pi}_i(D_i|X_i;\widehat{\alpha})}\cdot I_q(D_i)\cdot\Expect\left(Y_i(D_i)-\Psi^{-1} \left\{\cdot\right\}\biggr|I_q(D_i)=1, X_i\right)\right],\]
\[= \Expect_{X_i,D_i}\left[\pi_i(D_i|X_i;\alpha)\cdot\frac{1}{\widehat{\pi}_i(D_i|X_i;\widehat{\alpha})}\cdot I_q(D_i) \cdot \Expect\left(Y_i(D_i)-\Psi^{-1} \left\{ \cdot \right\}\biggr|I_q(D_i)=1,X_i\right)\right]\]
\[= \Expect_{D_i}\left[\Expect\left(Y_i(D_i)-\Psi^{-1} \left\{ \cdot \right\}\biggr|I_q(D_i)=1\right)\right]=0.\]
We can write this equation for each stratum $q$ in the form
\[\mu\left(\mathcal{D}_q\right)=\Expect_{j}\left[\Expect_i\left\{Y_i(d_{qj})\right\}\right]=\Expect_j \left[\Expect \left(\Psi^{-1} \left\{m\left(d_{qj},x_i; \beta \right)+\frac{\varphi_{q}} {\widehat{\pi}(d_{qj}|X_i;\widehat{\alpha}} \right\}  \right) \right].\]
Thus, averaging the mean values of the AOR model over $d_{qj} \in \mathcal{D}_q$ provides an unbiased estimate of $\mu\left(\mathcal{D}_q\right)$ regardless of whether $\Psi^{-1} \left\{m\left(D_i,X_i; \widehat{\beta} \right)\right\}$ is unbiased for $\Expect[Y_i|D_i,X_i]$.
\end{Thm}

The estimating equations of the DR model are implied by a number of standard estimators. For instance, Maximum Likelihood Estimation (MLE), Maximum Quasi-Likelihood (MQL), Restricted MLE (REML) for linear mixed models (LMMs), and Penalized Quasi-Likelihood (PQL) for generalized linear mixed models (GLMMs)

\subsection{Doubly robust dose-response estimation via Taylor series expansion}
The approach outlined above yields consistent estimates of the mean APO for each stratum. To obtain point estimates on the continuous dose-response function we apply a Taylor series expansion via a polynomial regression model. 

The true dose response function, $\mu(d)=\Expect\left[Y_i(d)\right]=f(d;\theta)$ is unknown, but we assume that it can be well approximated by an $Mth$ degree polynomial function
\[\Expect\left[Y_i(d)\right]\approx\sum_{m=0}^{M}\theta_m d^m. \] 
The mean APO for stratum $q$ relates to the polynomial dose-response in the following way
\[\Expect_{j}\left\{\Expect\left[Y_i(d_{qj})\right]\right\} \approx \sum_{m=0}^{M}\theta_m \Expect_{j}\left(d_{qj}^m \right),\]
with sample analogue 
\begin{equation}
\widehat{\mu}(\mathcal{D}_q) = \sum_{m=0}^{M}\theta_m \sum^{J}_{j=1} \frac{d_{qj}^m}{J},
\end{equation}
where $J$ is the number of units in stratum $q$. To evaluate this expression we need estimates of $\theta_m$. We obtain these by regressing the mean APO estimates for each stratum, $\widehat{\mu}(\mathcal{D}_q)$, on $M-1$ covariates $\sum_{j}^{J} d_{qj}^m/J$, $m=(1,...,M)$. Use of OLS will yield best linear unbiased estimates of $\theta_m$ under the Gauss Markov conditions. We use then use the OLS estimates of $\theta_m$ to recover a polynomial approximation to points on the dose response function via
\[\widehat{\mu}(d) = \sum_{m=0}^{M}\widehat{\theta}_m d^m.\]

Variances of the resulting dose-response point estimates can be obtained by applying a single bootstrap resampling scheme over the estimation of each model component (i.e. PGPS, AOR, polynomial regression).

\subsection{Comparison with existing doubly robust approaches}
To our knowledge there are currently only two papers that have constructed DR estimators using the GPS \citep[][construct a DR estimator for continuous treatments but in the context of marginal structural models and not for point exposure problems]{VanDerLaan/Robins:2003}.

\citet{Flores/Mitnik:2009} combine OR with inverse GPS weighting to obtain a DR weighted regression estimator for multi-valued, but not continuous treatments, in which doses take discrete values $d=1,...k$. They do this by weighting the OR model
\[\Expect[Y_i|D_i,X_i]=\sum_{j}^{k}\alpha_j\cdot 1_1 (D_i=j)+\beta^{\sf{T}}X_i\]
with weights $w_i=\sqrt{1/R_i}$. To calculate a weighted regression (WR) estimator at dose $d$ they use
\begin{equation}\label{FM}
\mu(d)_{WR1}= \frac{\sum_{i=1}^n \left(\widehat{\alpha}_d+\widehat{\beta}^{\sf{T}}x_i\right) \cdot w_i}{\sum_{i=1}^n w_i}.
\end{equation}
They do not provide a formal proof for the DR property but make two arguments for illustration. First, that weighting does not affect the consistency properties of the outcome regression; and second, that weighting by a correctly specified GPS gives a balanced sample analogous to adjusting for observed characteristics with experimental data.

\citet{Zhang/et/al:2012} construct a similar DR estimator for continuous treatments by combining and OR and a PS model through inverse probability weighting. Their weighted regression estimator is
\begin{equation}\label{Zhang}
\mu(d)_{WR2}=\frac{1}{N}\sum_{i=1}^{N}m(d,X_i,\widehat{\beta^*}),
\end{equation}
where $\widehat{\beta^*}$ is the estimate for the weighted regression model $m(D_i,X_i,\widehat{\beta^*})$ with weights $W(D_i)/R_i$ where $W(D_i)>0$ is chosen to stabilize the weights when the GPS takes very small values. In their application they chose $W(D_i)$ to be the marginal density of $D_i$ under a normal model.

The key difference between our augmented GPS approach, and the GPS weighting approaches, is that by including multiple inverse GPS covariates for different strata we are able to induce local rather than global bias correction. This follows from theorem 2 above where we show that the bias correction is specific to the mean APO for each stratum. If instead we were to employ a single weighting scheme, or include a single inverse GPS covariate, the bias correction would be tailored to a mean APO over the whole sample. For non-linear dose response-functions, local rather than global approximations allow to us to correct for misspecification due to an inaccurate assumed functional form in addition to effects of confounding. We illustrate these properties in the simulations presented in the next section.

\section{Simulations}\label{C6S4}
The model structure analysed in our simulations corresponds to the following data generating process.

A continuous treatment is assigned as a function of covariates $X_1$, $X_2$, and $Z$. The Gaussian dose-response function is quadratic in the treatment with confounding from $X_1$ and $X_2$.
\[X_1,X_2 \sim \mathcal{N}(\mu_{X_1}=4, \mu_{X_2}=8,\sigma^2_{X_1}=1,\sigma^2_{X_2}=2,\rho=-0.5)\]
\[Z \sim \mathcal{N}(10,4)\]
\[D\sim \mathcal{N}(2.0+0.5X_1+0.25X_2+Z, \sigma^{2}_{D})\]
\[Y\sim \mathcal{N}(1.0+4.0D-0.125D^2+0.5X_1+2X_2-0.5X_2^2,\sigma^{2}_{Y})\]
The correct OR model is
\[\Expect[Y|D,X_1,X_2]=\beta_0+\beta_{1}D+\beta_{2}D^2+\beta_{3}X_1+\beta_{4}X_2+\beta_{5}X_2^2,\]
and the correct PGPS model is
\[\widehat{\pi}_T^{-1}= \int_{D-\delta}^{D+\delta}\frac{1}{\sqrt{2\pi \widehat{\sigma}_D^2}}\exp\left(-\frac{1}{2\widehat{\sigma}_D^2}(t-\widehat{\mu}_D)^2 \right)dt,\]
where parameters $\mu_D$ and $\sigma_D^2$ are estimated from the correct regression model
\[\Expect[D|X_1,X_2]=\alpha_0 +\alpha_1 X_1+\alpha_2 X_2.\]

\subsection{Simulations for strata mean APO estimation}

The following models are tested:
\begin{itemize}
\item[1.] $\widehat{\mu}(\mathcal{D}_k)_{OR1}$ - based on the correctly specified OR model with predicted values averaged over $j=(1,...,J)$ treatment levels within each stratum
\[\widehat{\mu}_{OR1}(\mathcal{D}_k)=J^{-1}\sum_{j=1}^{J} \left[n^{-1} \sum_{i=1}^{n} \Psi^{-1} \left\{m(d_{kj}, X_i; \widehat{\xi})\right\}\right].\]
\item[2.] $\widehat{\mu}(\mathcal{D}_k)_{OR2}$ - same estimator as $\widehat{\mu}(\mathcal{D}_k)_{OR1}$ but based on an incorrectly specified OR model, with $X_1$ assumed as sole confounder.
\item[3.] $\widehat{\mu}(\mathcal{D}_k)_{DR1}$ - a DR estimator as described in equation (\ref{MUDR1}) above based on an incorrectly specified OR model ($X_1$ as sole confounder) augmented with correctly estimated inverse PGPS ($\widehat{\pi}_T^{-1}$) covariates for defined strata of the treatment.
\item[4.] $\widehat{\mu}(\mathcal{D}_k)_{DR2}$ - a DR estimator based on a correctly specified OR model augmented with incorrectly estimated inverse PGPS covariates ($\widehat{\pi}_F^{-1}$), with $X_1$ assumed as sole confounder, for defined strata of the treatment.
\item[5.] $\widehat{\mu}(\mathcal{D}_k)_{DR3}$ - a DR estimator based on an incorrectly specified OR model augmented with incorrectly estimated inverse PGPS covariates.
\item[6.] $\widehat{\mu}(\mathcal{D}_k)_{DR4}$ - a DR estimator based on an correctly specified OR model augmented with correctly estimated inverse PGPS covariates.
\end{itemize}

The simulations are based on 1000 runs on generated datasets of size 10,000. The mean of $d$ is 16 and the range approximately 1 to 30. Estimates are presented for the following treatment strata: (10.5,12.5], (12.5,14.5], (14.5,16.5], (16.5,18.5], (18.5,20.5]. To calculate the PGPSs we set $\delta=0.5$. Table \ref{sim1} shows our simulation results. The following results are reported: average estimates (Av Est), average estimated variances (Av Est Var), empirical variances (Emp Var) calculated from the bootstrap replications, mean squared error (MSE), and coverage based on normal bootstrap 95\% confidence intervals. 

The correctly specified OR model, OR1, provides unbiased estimates of points on the quadratic dose-response and consequently shows good results for all treatment intervals with relativity low variance and good coverage. The incorrectly specified OR model, OR2, fails to find the quadratic relationship between outcome and treatment levels, instead indicating a linear decreasing effect with poor coverage and large bias and MSE. Augmenting the incorrect OR model with the inverse PGPS covariates substantially corrects for bias from confounding and from functional misspecification of the treatment covariate. The DR1 estimator finds a quadratic relationship between outcome and treatment and performs well in terms of bias, MSE and coverage, although as we would expect, the variances are larger than for the correctly specified OR model (which has minimum variance amongst all consistent estimators). The DR2 model shows that the inclusion of addition irrelevant covariates in the AOR model does not induce bias, but as illustrated in theorem 1, it does increase the variance relative to the correct OR model. If the PGPS model and the OR model are wrong there is no DR property to be had and this is demonstrated in the DR3 results which indicate large bias and variance in estimation of the dose-response and poor coverage. Finally, if both the PS and the OR model are correctly specified the doubly robust estimator (DR4) provides unbiased estimates but with variance larger than would be achieved via a correct OR model alone. Thus, in summary, the finite sample properties of our DR approach indicate consistent estimation with relatively good performance in terms of efficiency. 

\begin{table}[htbp]
  \centering
 {\bf  \caption{Simulation results for Gaussian dose-response GLM with quadratic treatment effect. \label{sim1}}}
    \begin{tabular}{ll|ccccc}
 \hline
         && \multicolumn{4}{c}{treatment intervals} &  \\

        &  & (10.5,12.5] & (12.5,14.5] & (14.5,16.5] & (16.5,18.5] & (18.5,20.5] \\
 \hline
Truth &       & 15.426 & 17.176 & 17.926 & 17.676 & 16.426 \\
          &       &       &       &       &       &  \\
    OR1   & Av Est & 15.425 & 17.176 & 17.926 & 17.677 & 16.428 \\
          & Av Est Var & 0.012 & 0.010 & 0.010 & 0.010 & 0.011 \\
          & Emp Var & 0.011 & 0.010 & 0.010 & 0.010 & 0.011 \\
          & MSE   & 0.012 & 0.010 & 0.010 & 0.010 & 0.011 \\
          & Coverage & 95.400 & 95.200 & 94.800 & 94.500 & 95.200 \\
          &       &       &       &       &       &  \\
    OR2   & Av Est & 16.944 & 16.628 & 16.311 & 15.994 & 15.677 \\
          & Av Est Var & 0.025 & 0.015 & 0.010 & 0.012 & 0.019 \\
          & Emp Var & 0.016 & 0.007 & 0.004 & 0.008 & 0.019 \\
          & MSE   & 2.320 & 0.308 & 2.614 & 2.838 & 0.579 \\
          & Coverage & 0.000 & 0.417 & 0.000 & 0.000 & 0.000 \\
          &       &       &       &       &       &  \\
    DR1   & Av Est & 15.430 & 17.177 & 17.925 & 17.681 & 16.428 \\
          & Av Est Var & 0.072 & 0.046 & 0.037 & 0.038 & 0.050 \\
          & Emp Var & 0.072 & 0.046 & 0.037 & 0.038 & 0.050 \\
          & MSE   & 0.074 & 0.047 & 0.040 & 0.039 & 0.051 \\
          & Coverage & 94.500 & 94.400 & 95.500 & 94.400 & 94.900 \\
          &       &       &       &       &       &  \\
    DR2   & Av Est & 15.428 & 17.178 & 17.926 & 17.679 & 16.424 \\
          & Av Est Var & 0.024 & 0.017 & 0.016 & 0.016 & 0.019 \\
          & Emp Var & 0.024 & 0.018 & 0.014 & 0.016 & 0.019 \\
          & MSE   & 0.024 & 0.017 & 0.014 & 0.016 & 0.019 \\
          & Coverage & 95.833 & 94.792 & 94.792 & 96.250 & 94.375 \\
          &       &       &       &       &       &  \\
    DR3   & Av Est & 16.162 & 17.576 & 18.012 & 17.439 & 15.876 \\
          & Av Est Var & 0.072 & 0.045 & 0.037 & 0.040 & 0.059 \\
          & Emp Var & 0.080 & 0.043 & 0.037 & 0.042 & 0.062 \\
          & MSE   & 0.621 & 0.203 & 0.044 & 0.098 & 0.365 \\
          & Coverage & 25.833 & 51.875 & 91.042 & 77.708 & 38.125 \\
          &       &       &       &       &       &  \\
    DR4   & Av Est & 15.430 & 17.172 & 17.930 & 17.676 & 16.422 \\
          & Av Est Var & 0.024 & 0.018 & 0.016 & 0.016 & 0.020 \\
          & Emp Var & 0.025 & 0.017 & 0.016 & 0.016 & 0.020 \\
          & MSE   & 0.025 & 0.017 & 0.016 & 0.016 & 0.020 \\
          & Coverage & 94.100 & 93.000 & 95.000 & 94.800 & 94.800 \\
\hline
    \end{tabular}%
\end{table}%

\subsection{Simulation of existing doubly robust weighting approaches}

Table \ref{compsim} shows comparative results obtained by fitting the models of equation \ref{FM} (WR1) and equation \ref{Zhang} (WR2). The first thing worth noting is that when the OR model is correctly specified (DR2 models) both WR models provide consistent estimates, and this is because theorem 1 above holds. Note, however, that the variance of the WR1-DR2 estimates is relatively large. When the OR model is incorrectly specified however, the WR estimators based on a global weighting scheme produce biased estimates of the mean APO for each stratum. Due to the inclusion of dummies in WR1-DR1 model, a quadratic shape is implied in the strata estimates, but the global bias correction induced by weighting does not adequately address local misspecification. With no dummies, as in WR2-DR1, the model fails to find the quadratic shape of the dose response and again mean APO estimates are biased. Note that excluding dummies from the WR1 approach we get estimates that are approximately equivalent to those of WR2.

\begin{table}[htbp]
  \centering
 {\bf  \caption{Comparative results based on GPS weighting estimators. \label{compsim}}}
    \begin{tabular}{ll|ccccc}
 \hline
         && \multicolumn{4}{c}{treatment intervals} &  \\

        &  & (10.5,12.5] & (12.5,14.5] & (14.5,16.5] & (16.5,18.5] & (18.5,20.5] \\
 \hline
Truth &       & 15.426 & 17.176 & 17.926 & 17.676 & 16.426 \\
&&											      &       &       &       &  \\
    WR1-DR1 & Av Est & 15.854 & 17.409 & 17.965 & 17.561 & 16.160 \\
          & Av Est Var & 0.068 & 0.044 & 0.037 & 0.039 & 0.053 \\
          & Emp Var & 0.065 & 0.044 & 0.036 & 0.040 & 0.056 \\
          & MSE   & 0.248 & 0.098 & 0.037 & 0.054 & 0.127 \\
          & Coverage & 63.750 & 82.083 & 95.625 & 91.458 & 79.583 \\
		  &       &       &       &       &       &  \\
    WR1-DR2 & Av Est & 15.420 & 17.171 & 17.921 & 17.665 & 16.420 \\
          & Av Est Var & 0.023 & 0.017 & 0.016 & 0.016 & 0.020 \\
          & Emp Var & 0.022 & 0.017 & 0.015 & 0.017 & 0.019 \\
          & MSE   & 0.022 & 0.017 & 0.015 & 0.017 & 0.019 \\
          & Coverage & 94.792 & 95.208 & 95.625 & 95.625 & 95.417 \\
		  &       &       &       &       &       &  \\
    WR2-DR1 & Av Est & 16.700 & 16.525 & 16.351 & 16.176 & 16.002 \\
          & Av Est Var & 0.012 & 0.009 & 0.009 & 0.011 & 0.015 \\
          & Emp Var & 0.012 & 0.010 & 0.009 & 0.011 & 0.015 \\
          & MSE   & 1.635 & 0.433 & 2.491 & 2.261 & 0.196 \\
          & Coverage & 0.000 & 0.000 & 0.000 & 0.000 & 6.875 \\
		  &       &       &       &       &       &  \\
    WR2-DR2 & Av Est & 15.420 & 17.184 & 17.941 & 17.692 & 16.437 \\
          & Av Est Var & 0.012 & 0.010 & 0.010 & 0.010 & 0.011 \\
          & Emp Var & 0.012 & 0.011 & 0.011 & 0.011 & 0.012 \\
          & MSE   & 0.012 & 0.011 & 0.011 & 0.012 & 0.012 \\
          & Coverage & 94.700 & 95.100 & 95.000 & 95.300 & 95.600 \\
\hline
    \end{tabular}%
\end{table}%

\subsection{Testing for joint significance of the AOR inverse GPS covariates}

In contrast to doubly robust weighting approaches, use of an AOR model allows tests to be performed on the inclusion of inverse GPS covariates for bias correction. In table \ref{wald} below we present results for Wald Chi-squared tests for the joint significance of the inverse GPS covariates in our DR simulations.

\begin{table}[htbp]
  \centering
 {\bf  \caption{Simulation of Wald Chi-squared tests for joint significance of inverse GPS covariates, rejections rates for the null: $H_0: \varphi=0$  (nominal level is 95\%) \label{wald}}} 
    \begin{tabular}{lc}
   \hline
          & $H_0$ rejection rate \\
\hline
    DR1   & 100.0\% \\
    DR2   & 4.7\% \\
    DR3   & 100.0\% \\
    DR4   & 4.8\% \\
\hline
    \end{tabular}%
\end{table}%

When the OR model is correctly specified (i.e. DR2 \& DR4) the rejection rate of the null hypothesis ($\varphi=0$) is approximately 5\% given either a correctly or incorrectly specified GPS model. The size of the test is thus as expected. When the OR model is incorrectly specified (i.e. DR1 \& DR3) the rejection rate of the null is 100\%, regardless of whether the GPS model is correctly specified or not. Thus, while the test shown here cannot provide inference on how well the GPS covariates have helped to reduce bias from confounding, it can help indicate when their inclusion offers little beyond an OR specification. This is important, because as shown above, inclusion of irrelevant covariates in the AOR model can reduce efficiency.     

\subsection{Simulations for dose-response estimation}

We now illustrate how Taylor expansion via polynomial regression can be used to recover an approximation to the dose-response curve using our mean APO estimates. For the same set up described above we obtain 1000 simulated set of results for different treatment strata, sliding the partition of the treatment to increase coverage of mean APOs. We then take the mean APO estimates for each strata as a response variable and fit a second order polynomial model to the mean dose as described in section 3 above. The parameters of the polynomial model are then used to recover an approximation to the dose response. 

\begin{figure}[htp]
\centering {\includegraphics[scale=0.5]{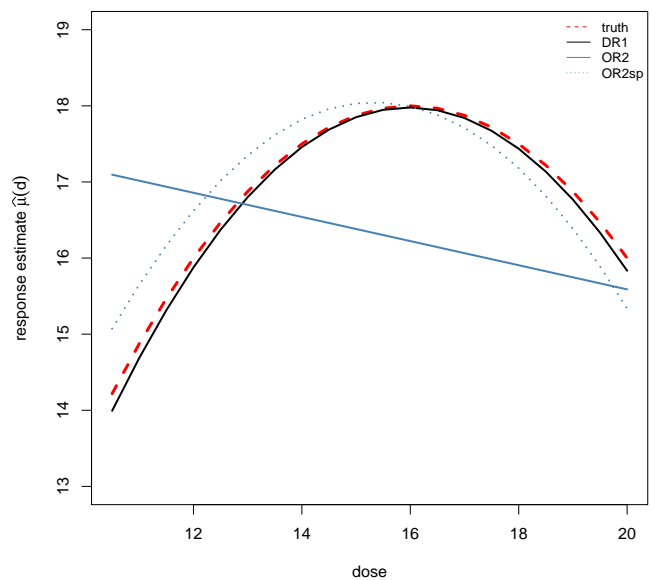} }
   {\bf \caption{Estimated dose-response functions obtained from fitting second order polynomial models to the DR1 and OR2 mean APO estimates. \label{C6F1}}}
\end{figure}

The results are illustrated graphically in figure \ref{C6F1} which shows the resulting dose-response curve obtained for DR1 and OR2 estimates relative to the truth. The figure also shows the mean dose-response curve obtained when the OR2 model is fitted using a Generalised Additive Model (GAM) (OR2sp), in which case model misspecification arises from confounding rather than functional misspecification.

The nonparametric fit to the DR1 estimates recovers a good approximation to the true dose-response curve. The OR2 model finds an incorrect shape for the dose-response because it is based on biased estimates of the mean APO for each stratum as discussed above. Fitting the OR2 model semiparametrically, with a smooth term for treatment $d$, gives a quadratic dose-response curve but which is still biased due to the effect of confounding from omitted covariates. This comparison thus demonstrates the ability of the DR model to correct for misspecification bias arising from both confounding and an incorrect functional form.

\section{Conclusions}
In this paper we have constructed a doubly robust approach for estimation of dose-response functions for continuous treatments. The model is doubly robust in the sense that consistent estimates of average potential outcomes can be obtained for defined strata of the treatment under misspecification of either the outcome regression or propensity score models. This is achieved by inducing local bias correction for misspecification of the outcome regression model by augmenting it with inverse propensity score covariates corresponding to defined strata of the treatment. A polynomial regression is then fitted to these point estimates to derive a Taylor approximation to the continuous dose-response function. Standard errors are derived by bootstrapping over all model components.

We have shown that our approach can provide a good approximation to linear or nonlinear dose-response curves and is robust to problems of confounding or functional form misspecification. Furthermore, simulations show that the efficiency performance in finite samples is good relative to minimum variance consistent estimators. There are several issues in the paper that require further attention in future research. In particular, the choice of treatment strata for point estimation could be given a more formal basis. Similarly, the choice of delta for estimation of propensity score probabilities should be examined in greater depth than has been possible here.

\bibliographystyle{chicago}
\bibliography{D:/Bibtexfiles/causal}

\end{document}